\documentclass[10 pt,leqno]{amsart}
\usepackage {amsfonts}
\usepackage{amsthm}\usepackage{amssymb}
\usepackage{latexsym}
\usepackage{amsmath}
\usepackage{mathrsfs}
\usepackage{pifont}
\theoremstyle{plain}
\newtheorem{tw}{Theorem}[section]

\newtheorem {lem} [tw]{Lemma}
\newtheorem {prop}[tw] {Proposition}

\newtheorem{cor}[tw]{Corollary}

\theoremstyle{definition}

\newcommand{\Ex} {\Bbb E}

\newcommand{\N}{\Bbb N}

\newcommand{\Z}{\Bbb Z}

\newcommand{\Ad}{\textrm{Ad}}

\newcommand{\Ll} {\mathcal L}

\newcommand{\Tt} {\mathcal T}
\newcommand{\Oo} {\mathcal O}
\newcommand{\Ff} {\mathcal F}
\newcommand{\Kk} {\mathcal K}

\newcommand{\id}{I}

\newcommand{\semidir}{{\Bbb o}}
\newcommand{\TE}{\Tt_E}
\newcommand{\Oe}{\Oo_E}
\newcommand{\OX}{\Oo_X}
\newcommand{\GE}{\Gamma_E}
\newcommand{\GXH}{\hat{\Gamma}_X}
\newcommand{\subs}{\subseteq}

\newcommand{\la}{\langle}
\newcommand{\ra}{\rangle}
\newcommand{\lla}{\left\langle}
\newcommand{\rra}{\right\rangle}
\newcommand{\eps}{\varepsilon}
\newcommand{\ot}{\otimes}

\newcommand{\wt}{\widetilde}

\subjclass[2000]{ Primary 46L05, Secondary 46B28}

\begin{document}
\author{Adam Skalski}
\author{Joachim Zacharias}

\footnote{\emph{Permanent address of the first named author:}
Department of Mathematics, University of \L\'{o}d\'{z}, ul.Banacha 22, 90-238 \L\'{o}d\'{z}, Poland.}
\address{Mathematics and Statistics, Lancaster University, Lancaster, LA1 4YF}
\email{a.skalski@lancaster.ac.uk}
\address{School of Mathematical Sciences,  University of Nottingham,Nottingham, NG7 2RD}
\email{joachim.zacharias@nottingham.ac.uk  }

\keywords{$C^*$-algebra, approximation properties, Hilbert bimodule, crossed product, Pimsner algebra, completely positive lift}

\title{\bf On approximation properties of Pimsner algebras and crossed products by Hilbert bimodules}

\maketitle

\begin{abstract}
\noindent
Let $X$ be a Hilbert bimodule over a $C^*$-algebra $A$ and $\Oo_X= A \rtimes_X \Z$.
Using a finite section method we construct a sequence of completely positive contractions
factoring through matrix algebras over $A$ which act on $s_{\xi} s_{\eta}^*$ as Schur multipliers converging to the identity. This shows immediately that for a finitely generated $X$ the algebra $\Oo_X$ inherits any standard approximation property such as nuclearity, exactness, CBAP or OAP from $A$.
We generalise this to certain general Pimsner algebras by proving semi-splitness of the Toeplitz extension under certain conditions and discuss some examples.
\end{abstract}

\vspace*{1cm}

\section*{Introduction}
In recent years Pimsner algebras (\cite{pi}) have been studied intensively since they form a rich class of algebras, containing  e.g.\;crossed products by $\Z$ and $\N$ and various Cuntz-type algebras and combine in a very flexible
way standard constructions of forming new $C^*$-algebras from old ones.
Also approximation properties of Pimsner algebras have been considered:
in \cite{ds} it was shown that an extended Pimsner algebra is exact iff its coefficient algebra is exact, the same result holds true for nuclearity (\cite{ge}) and the completely bounded approximation property
(cf.\;\cite{dsm}, where it is also shown that the Haagerup constant of both algebras are the same). All known proofs of
approximation properties of Pimsner algebras use a crossed product representation of a certain dilated algebra and refer, in effect,
to known results for crossed products by a single automorphism, in particular
no explicit approximations are given.

In this paper we construct explicit approximations using an idea from the theory of Toepliz operators, namely the finite section method which works best for Pimsner algebras over Hilbert bimodules (i.e.\;crossed products).
Our methods also allow to prove slightly stronger results than in \cite{ds} and \cite{dsm} in certain situations, for instance for the OAP. (Although these could probably also be obtained using the crossed reprsentations.)

In particular
we obtain a sequence of completely positive contractions
factoring through matrix algebras over $A$ which act on $s_{\xi} s_{\eta}^*$ as Schur multipliers converging to the identity.

In a broad sense the finite section method studies Toeplitz operators on
say $\ell^2(\N)$ in terms of their compressions to the subspaces $\textrm{span}\{e_1 , \ldots ,e_n\}$,
where $(e_i)_{i=1}^{\infty}$ denotes a standard orthonormal basis of $l^2$. Asymptotically every Toeplitz operator
can be reconstructed from its compressions. It turns out that in general if $E$ is a $C^*$-correspondence, each element
in the Toeplitz algebra $\Tt_E$ acting on the Fock module $\Gamma_E$ can asymptotically be reconstructed from its finite
sections in a completely positive way at least modulo $\Kk(\Gamma_E)$. This means that we can find an approximation of the
quotient map from the Toeplitz-Pimsner algebra of $E$ to the Cuntz-Pimsner algebra of $E$. Since it is easy to see that for
a finitely generated $E$
a finite section algebra has the same approximation properties as the coefficient algebra it follows immediately that both,
the Cuntz-Pimsner and Toeplitz-Pimsner algebras, have all the approximation properties of the coefficient algebra provided there exists a completely positive contactive lift of the quotient map.
Pimsner constructed such a lift in his original paper in the case where $E$
is a Hilbert bimodule. We extend this result by showing that such a lift exists
when $E$ is only a $C^*$-correspondence proviede there is a conditional expectation of the algebra $\mathcal F_E$ onto $A$.

The structure of the paper is as follows: in the first section the finite section method for Pimsner algebras is described.
 Section 2 contains the discussion of the existence of lifts of the exact sequence
$0 \to  \Kk (\Gamma_E) \to \TE \to \Oe \to 0$. Finally in Section 3 various approximation properties of $C^*$-algebras
are recalled and the contents of the first two sections are applied to prove the main results of the paper.

\section{The Finite Section Method for Pimsner Algebras}            \label{finitesection}

Throughout the paper $A$ will denote a unital $C^*$-algebra and $E$ a \emph{$C^*$-correspondence} over $A$,
i.e.\;a Hilbert $A$-module together with a nondegenerate faithful left action $\phi: A \to \Ll_A(E)$.
A particular example is a \emph{Hilbert bimodule} (\cite{bcm}) i.e.\;a bimodule with two different scalar products
$\mbox{}_A\la \cdot , \cdot \ra$ and $\la \cdot , \cdot \ra_A$, such that for all $\xi, \eta, \zeta \in E$
$$\xi \la \eta , \zeta \ra_A= \mbox{}_A\la \xi , \eta \ra \zeta.$$
A Hilbert bimodule satisfies the same compatibility conditions as an equivalence bimodule, but the corresponding left and right Hilbert modules are not required to be full (cf.\;\cite{bcm}). Hilbert bimodules will usually be denoted by $X$.

In his seminal paper \cite{pi} Pimsner defined two $C^*$-algebras associated to $E$. The Toeplitz-Pimsner algebra $\Tt_E$ is the concrete $C^*$-algebra generated by Toeplitz operators on the Fock module $\Gamma_E= \bigoplus_{n\geq 0} E^n$ defined by the formula $t_{\xi} \eta = \xi \ot \eta$, where $\xi \in E$ and $\eta \in \Gamma_E$. Here and subsequently $E^0:=A$ and $E^n:= E^{\ot n}$ denotes the $n$-fold internal tensor product (see \cite{Lance}) for $n \geq 1$. The Cuntz-Pimsner algebra $\Oo_E$  is the quotient $\Tt_E / \Kk(\Gamma_E)$, giving the following exact sequence:
\begin{equation} \label{exact}
0 \to \Kk (\Gamma_E) \to \TE \to \Oe \to 0.\end{equation}
We denote the generators of $\Oe$, arising as images of $t_{\xi}$ ($\xi \in E$) in the quotient map, by $s_\xi$; this notation is also extended in a natural way for $\xi \in E^r$, where $r\in \N$. We write $|\xi| = r$ if $\xi \in E^r$.

In case of a Hilbert bimodule $X$ one can also consider a two-sided Fock module $\hat{\Gamma}_X = \bigoplus_{n\in \Z} X^n$, where $X^{-1}=X^*, \; X^{-2} = (X^2)^* , \ldots$ It may be shown that left creation operators on $\hat{\Gamma}_X$ generate the $C^*$-algebra $\Oo_X$ (see Section \ref{liftings}). It turns out that $\Oo_X$ may also  be regarded as a crossed product of $A$ by a Hilbert bimodule {\cite{aee}} and the two-sided Fock representation may be considered as the regular representation of $\Oo_X= A \semidir_X \Z$. In Section \ref{liftings} it is explained how a simple compression given by the projection from $\hat{\Gamma}_X$ onto $\Gamma_X$ provides a ucp (ucp here and in what follows stands for unital completely positive) lift to the exact sequence \eqref{exact}, in case where $E$ is replaced by a Hilbert bimodule $X$.

The dense $*$-subalgebra of $\TE$ generated by $A$ and $t_{\xi}$ ($\xi \in E$) is spanned by the canonical elements $t_{\mu} t_{\nu}^*$, where $\mu $ and $\nu$ lie in arbitrary powers of $E$ (including $E^0$). Using the `rank one' operators $e_{\mu , \nu} \in \Kk(E^{|\nu|} , E^{|\mu|})$ defined by the formula
\[ e_{\mu , \nu} (\xi) = \mu \la \nu , \xi \ra, \;\;\;\xi \in E^{|\nu|} \]
(in the literature often denoted by $\Theta_{\mu , \nu}$), we obtain the equality
$$t_{\mu} t_{\nu}^* = \sum_{k=0}^{\infty} e_{\mu , \nu} \otimes \id_{E^k},$$
where convergence of the infinite sum is understood in the strict topology of $\Ll_A(\GE)$.
In case of a Hilbert bimodule $X$ we obtain the formula
$$s_{\mu} s_{\nu}^* = \sum_{k=-\infty}^{\infty} e_{\mu , \nu} \otimes \id_{X^k},$$
for the elements $s_{\mu} s_{\nu}^* \in \Oo_X \subset \Ll_A(\hat{\Gamma}_X)$.

For $N \in \N$ let $\GE^{(N)} = A \oplus E \oplus \ldots \oplus E^N \subseteq \GE$ and let
$P_N \in \Ll_A( \GE)$ be the corresponding projection. Define the completely positive compression $\phi_N : \TE \to \Ll_A ( \GE^{(N)})$ by
$$\phi_N (x) = P_N x P_N, \;\; x \in \TE.$$
Every element  $x \in \Ll_A (\GE^{(N)})$ has a unique decomposition $x= \sum_{i,j=0}^{N} x_{i,j}$, where $x_{i,j} \in \Ll_A(E^j,E^i)$. Since $x_{i,j} \otimes \id_{E^k} \in \Ll_A(E^{j+k} , E^{i+k} )$ we may define $x \otimes \id_{E^k} = \sum_{i,j} x_{i,j} \otimes \id_{E^k}$ and moreover consider the map $ \Psi_N: \Ll_A (\GE^{(N)}) \to  \Ll_A (\GE)$ given by the formula ($ x \in \Ll_A (\GE^{(N)})$)
$$\Psi_N (x) = (N+1)^{-1} \sum_{k=0}^{\infty} x \otimes  \id_{E^k}.$$
The map $\Psi_N$ is completely positive and since
\begin{eqnarray*}
\Psi_N (1)&=&\Psi_N( 1_A \oplus \id_E \oplus \ldots \oplus \id_{E^N} ) \\
&=&(N+1)^{-1} \left( \sum_{k=0}^{\infty} \id_{E^k} +\sum_{k=1}^{\infty} \id_{E^k} + \ldots + \sum_{k=N}^{\infty} \id_{E^k} \right) \\
&=&(N+1)^{-1} \left( \sum_{k=1}^{N+1} k \id_{E^{k-1}} \right) \otimes \id_{\GE}\end{eqnarray*}
it is also a contraction.

\begin{lem}       \label{quotappr}
For each $N \in \N$ the map $\Psi_N \circ \phi_N :\Ll_A (\GE) \to \Ll_A (\GE)$ is a completely positive contraction mapping $\TE$ into $\TE$. Let $p: \Tt_E \to \Oo_E$ be the canonical quotient map, $V_N= p \circ \Psi_N \circ \phi_N$ and $\mu \in E^{|\mu|}$, $\nu \in E^{|\nu|}$. Then
$$
V_N(t_{\mu}t_{\nu}^*)=
\left\{\begin{array}{ccc}\frac{\min (N-|\mu| , N-|\nu |)}{N+1}  s_{\mu}s_{\nu}^* & \textrm{if} & | \mu| , |\nu| \leq N , \\
0 &  & \textrm{otherwise.}
\end{array}
\right.
$$
Moreover, $V_N (x) \to p (x) $ as $N \to \infty$ for each $x \in \TE$.
\end{lem}

\begin{proof}
Essentially, we only need to verify the claimed identity. Since $P_N t_{\mu} =0$ if $| \mu | > N$ it is clear that $V_N (t_{\mu}t_{\nu}^*)=0$ whenever $|\mu|>N$ or $|\nu | >N$. In case $|\mu| , | \nu | \leq N$
$$P_N t_{\mu}t_{\nu}^*P_N = P_N \left( \sum_{k=0}^{\infty} e_{ \mu , \nu} \otimes \id_{E^k} \right)P_N = \sum_{k=0}^{\min (N-|\mu| , N -|\nu|)} e_{ \mu , \nu} \otimes \id_{E^k}.$$
Since $\Psi_N \big( e_{ \mu , \nu} \otimes \id_{E^k} \big) = (N+1)^{-1}\Big( \sum_{l=k}^{\infty} e_{ \mu , \nu} \otimes \id_{E^l} \Big)$ we obtain
$$\Psi_N \circ \phi_N(t_{\mu}t_{\nu}^*) = \tfrac{\min (N-|\mu| , N-|\nu |)}{N+1} t_{\mu}t_{\nu}^* \; - \;(N+1)^{-1} \left( \sum_{l=1}^{\min( N - |\mu| , N- |\nu|)} \sum_{l=0}^k e_{\mu ,\nu} \otimes \id_{E^l} \right)$$
The second term lies in $\Kk(\GE)$, which allows to conclude the proof.
\end{proof}

Thus the finite section method yields an approximation of the quotient map $p$, where each of the approximating maps $V_N$ ($N \in \N$) acts like a Schur multiplier. The approximation problem for $\TE$ and $\Oe$ is thereby reduced to finding a completely positive lifting of $p$. We will address this problem in the following section. As was mentioned before, in case $E$ is a Hilbert bimodule, Pimsner has constructed such a lift in a somewhat abstract way. However, our method works directly in this case and no lifting procedure is needed. Let $X$ be a Hilbert bimodule over $A$. In this case $\OX$ acts on $\GXH$ by $s_{\xi} \eta = \xi \ot \eta$ and we may modify the definitions of approximating maps in an obvious way. Let ($N \in \N$) $\GXH^{(N)}=A \oplus X \oplus \ldots \oplus X^N$ with the corresponding projection again denoted by $P_N$, let $\hat{\phi}_N (x) = P_N x P_N$ for $x \in \Oo_X$ and define $\hat{\Psi}_N:  \Ll_A (\GXH^{(N)}) \to  \Ll_A (\GXH)$ by the formula
$$\hat{\Psi}_N (x)=(N+1)^{-1} \sum_{k=-\infty}^{\infty} x \otimes  \id_{X^k},\;\; x \in \Ll_A (\GXH^{(N)}).$$
 Since $\id_{X^l} \ot \id_{X^k} = \id_{X^{l+k}}$ for $k,l \in \Z$ we have $\hat{\Psi}_N (\id_{\GXH^{(N)}})=\id_{\hat{\Gamma}_X}$ i.e.\ $\hat{\Psi}_N $ is unital. Then the following is an easy modification of Lemma \ref{quotappr}.

\begin{lem}
 Let $N \in \N$. The map $W_N= \hat{\Psi}_N \circ \hat{\phi}_N: \Ll_A (\GXH) \to \Ll_A (\GXH)$is a completely positive contraction mapping $\Oo_X$ into $\Oo_X$ and, for $\mu \in X^{|\mu|}$, $\nu \in X^{|\nu|}$,
$$
W_N( s_{\mu}s_{\nu}^*)=
\left\{\begin{array}{ccc}\frac{\min (N-|\mu| , N-|\nu |)}{N+1}  s_{\mu}s_{\nu}^* & \textrm{if} & | \mu| , |\nu| \leq N , \\
0 &  & \textrm{otherwise.}
\end{array}
\right.
$$
Moreover, $W_N (x) \to x$ as $N \to \infty$ for $x \in \OX$.
\end{lem}
\begin{proof}
Since $\sum_{k=-\infty}^{\infty} e_{\mu,\nu} \otimes \id_{X^l} \ot \id_{X^k} =\sum_{k=-\infty}^{\infty} e_{\mu,\nu} \otimes \id_{X^k}$ we have no perturbation term and the formula follows.
\end{proof}

Here each $W_N$ acts precisely as a Schur type multiplier on the canonical generators.

\section{Liftings} \label{liftings}

\noindent
In this section we discuss (ucp) liftings of the exact sequence
$$0 \to  \Kk (\Gamma_E) \to \TE \to \Oe \to 0.$$
It does not seem to be known whether such liftings always exist. Below we consider some special cases. As already pointed out in case $E$ is a Hilbert bimodule $X$ (over $A$) such a lift has been constructed in \cite{pi} as follows. Let again $s_{\xi} = t_{\xi} + \Kk (\Gamma_X)$ ($\xi \in X$) denote the generators of ${\mathcal O}_X = T_X / \Kk (\Gamma_X)$ and let $\bar{s}_{\xi} $ denote the left creation operators on the two-sided Fock module $\hat{\Gamma}_X$, generating the $C^*$-algebra $\bar{\Oo}_X= C^*(\bar{s}_{\xi} , a \mid \xi \in X , a \in A)$. Then by universality of $\Oo_X$ there exists a $*$-homomorphism $\Oo_X \to \bar{\Oo}_X$ sending $s_{\mu} s_{\nu}^*$ to $\bar{s}_{\mu} \bar{s}_{\nu}^*$.
A completely positive map $\bar{\Oo}_X \to \Tt_X$, $x \to PxP$ (where $P$ is the projection of $\GXH$ onto $\Gamma_X$) sending $\bar{s}_{\mu} \bar{s}_{\nu}^*$ to $t_{\mu} t_{\nu}^*$ and finally the quotient map $\Tt_X \to \Oo_X$ sending $t_{\mu} t_{\nu}^*$ to $s_{\mu} s_{\nu}^*$ for all $r, s \in \N,$ $ \mu \in X^r, \nu \in X^s$ complete the diagram:
$$\Oo_X \to \bar{\Oo}_X \to \Tt_X \to \Oo_X.$$
Since the composition of all three maps yields  the identity map and the first map is surjective, we see that $\Oo_X \cong \bar{\Oo}_X$ (a fact that has already been used in Section \ref{finitesection}). Composition of the first two maps gives a desired ucp lift for the sequence \eqref{exact}.

Given any finitely generated $C^*$-correspondence $E$ let $\Ff_E = \lim_{n \to \infty} \Kk (E^n)$ be an inductive limit $C^*$-algebra, given by the usual embeddings $\Kk(E^j) \hookrightarrow \Kk(E^{j+1})$ defined by $T \to T \ot 1$. Note that $A$ may be identified with a $C^*$-subalgebra of $\Ff_E$, and the internal tensor product (over $A$)  $E \otimes \Ff_E$ is a Hilbert bimodule over $\Ff_E$ under the scalar products ($\xi_1, \xi_2 \in E, b_1, b_2 \in \Ff_E$)
$$
\mbox{}_{\Ff_E} \langle \xi_1 \ot b_1 , \xi_2 \ot b_2 \rangle = e_{\xi_1 , \xi_2 } \otimes b_1 b_2^*\;\;\;\; \text{ and } \;\;\;\;\; \langle \xi_1 \ot b_1 , \xi_2 \ot b_2 \rangle_{\Ff_E} =  b_1^* \langle \xi_1 , \xi_2 \rangle  b_2 .$$
By \cite{pi} $\Oo_{E \ot \Ff_E}  \cong \Oo_E$ always holds, but for Toeplitz-Pimsner algebras in general we only have an inclusion $\Tt_E \subs \Tt_{E \ot \Ff_E}$.

Now let $B$ be a $C^*$-algebra such that $A \subs B$ and let $\varepsilon : B \to A$ be a conditional expectation.

\begin{lem} \label{condexp}
For any Hilbert $A$-module $E$ the map $\xi \ot b \to \xi \eps (b)$ extends to an $A$-linear contraction $\bar{\eps}: E \ot B \to E$ (note $E \ot B$ is a $B$-Hilbert module).
\end{lem}
\begin{proof}
Given $ n\in \N$, $b_1, \ldots,b_n \in B$ and $\xi_1 ,\ldots ,\xi_n \in E$, we have
\begin{align*}\Big\langle  \bar{\eps}(\sum_{i=1}^n \xi_i \ot b_i ) ,\bar{\eps}( \sum_{j=1}^n \xi_j \ot b_j ) \Big\rangle_E& = \sum_{i,j=1}^n \eps (b_i)^* \langle \xi_i , \xi_j \rangle \eps (b_j) \\
&= [\eps (b_1)^* \dots \eps (b_n)^*] [\langle \xi_i , \xi_j \rangle ] [\eps (b_1) \dots \eps (b_n)]^t.\end{align*}
Notice that by the Schwarz inequality for completely positive maps and properties of conditional expectations $\eps (b^*a^*ab) \geq \eps(b)^* a^*a \eps(b)$ for $a \in A$ and $b \in B$. A similar statement holds for $\eps^{(n)}= \eps \ot \textup{id}_{M_n}$ since $\eps^{(n)}$ is also a conditional expectation. Applying the latter to the matrix  $[\langle \xi_i , \xi_j \rangle]_{i,j=1}^n \in M_n(A)$ yields
\begin{eqnarray*}\sum_{i,j=1}^n \eps (b_i)^* \langle \xi_i , \xi_j \rangle \eps (b_j)& \leq &\sum_{i,j=1}^n \eps \left( b_i^* \langle \xi_i , \xi_j \rangle b_j \right) \\
&=&\eps \left( \Big\langle \sum_{i=1}^n \xi_i \ot b_i , \sum_{j=1}^n \xi_j \ot b_j \Big\rangle \right)\end{eqnarray*}
and so
$$\Big\| \Big\langle  \sum_{i=1}^n \xi_i \eps (b_i) , \sum_{j=1}^n \xi_j \eps (b_j) \Big\rangle \Big\| \leq \Big\|\Big\langle \sum_{i=1}^n \xi_i \ot b_i , \sum_{j=1}^n \xi_j \ot b_j \Big\rangle \Big\|,$$
i.e.\ the prescription $\xi \ot b \mapsto \xi \eps (b)$ indeed extends (by linearity and continuity) to a contraction $\bar{\eps}$.
\end{proof}

It is well known that whenever $E_1,E_2$ are Hilbert modules  ($E_1$ over a $C^*$-algebra $C$) and the internal tensor product $E_1 \ot E_2$ is defined, the tensor product $\phi \ot \id_{E_2}$ is a well defined operator for every completely positive map in $\Ll_C(E_1)$ (\cite{Lance}). In the above lemma we need the conditional expectation property, as we define a map via tensoring by identity on the left: $\id_E \ot \eps$. For the further use note the following facts: first, we have the equality:
\begin{equation} \label{first}
\langle \xi, \bar{\eps} \zeta \ra_{E} = \eps \left( \la \xi \ot 1, \zeta \ra_{E \ot B} \right), \;\;\; \xi \in E, \zeta \in E \ot B.\end{equation}
Secondly, for all $n \in \N$ the map $ \bar{\eps}_n = \bar{\eps} \oplus \cdots \oplus \bar{\eps}: (E \ot B)^{\oplus n} \to E ^{\oplus n} $ via the natural isomorphism $(E \ot B)^{\oplus n} \cong E \ot B^{\oplus n} $ (and viewing $ E ^{\oplus n}$ as an $A^{\oplus n} $ module) is identified with the map $\bar{\eps'}$ introduced as was done in Lemma \ref{condexp} but this time for the conditional expectation $\eps'=   \eps \oplus \cdots \oplus \eps :B^{\oplus n}  \to A^{\oplus n} $. Therefore formula \eqref{first} has the `matricial' counterpart
\begin{equation} \label{second}
\langle \xi, \bar{\eps}_n \zeta \ra_{E} =\eps' \left( \la \xi \ot 1, \zeta \ra_{(E \ot B)^{\oplus n} } \right), \;\;\;\xi \in E^n, \zeta \in (E \ot B)^{\oplus n} .\end{equation}

The above allow us to deduce the following lemma.

\begin{lem} \label{hatexp}
The map $\hat{\eps}: \Ll_B(E \ot B) \to \Ll_A (E)$ defined by the formula
\[ \hat{\eps} (T) (\xi) = \bar{\eps} \left(T (\xi \ot 1) \right), \;\;\; T \in \Ll_B (E \ot B), \xi \in E,\]
is unital and completely positive. For all $\xi ,\eta \in E$ and $b,c \in B$,
\[ \hat{\eps}(e_{\xi\ot b,\eta \ot c} ) =e_{\xi \eps(b c^*),\eta } = e_{\xi ,\eta \eps(c b^*)} .\]
In particular $\hat{\eps}(\Kk(E \ot B))= \Kk(E)$.
\end{lem}
\begin{proof}
The fact that $\hat{\eps}$ is well defined, unital and bounded follows  directly from $A$-linearity of $\bar{\eps}$ and elementary checks. Complete positivity may be also checked directly: let $n \in \N$, $[T_{ij}]_{i,j=1}^n \in \Ll_B (E \ot B) \ot M_n \cong  \Ll_B ((E \ot B)^{\oplus n})$ be positive and $\xi_1, \ldots, \xi_n\in E$. Then formula \eqref{second} implies the following (see the notations introduced before the lemma):
\begin{align*}  \lla \begin{bmatrix} \xi_1 \\ \vdots \\ \xi_n \end{bmatrix},\hat{\eps}^{(n)} \left([T_{ij}]_{i,j=1}^n \right)\begin{bmatrix} \xi_1 \\ \vdots \\ \xi_n \end{bmatrix}  \rra =&\lla \begin{bmatrix} \xi_1 \\ \vdots \\ \xi_n \end{bmatrix},\bar{\eps}_n \left( [T_{ij}]_{i,j=1}^n \begin{bmatrix} \xi_1 \ot 1 \\ \vdots \\ \xi_n \ot 1 \end{bmatrix} \right) \rra \\
=&\, \eps' \left( \lla \begin{bmatrix} \xi_1 \ot 1 \\ \vdots \\ \xi_n \ot 1 \end{bmatrix}, [T_{ij}]_{i,j=1}^n \begin{bmatrix} \xi_1 \ot 1 \\ \vdots \\ \xi_n \ot 1 \end{bmatrix} \rra \right) \geq 0,\end{align*}
and $\hat{\eps}$ is completely positive. Moreover for $\xi ,\eta, \zeta \in E$ and $b,c \in B$,
\begin{eqnarray*}\hat{\eps}(e_{\xi \ot b, \eta \ot c})( \zeta)&= &\bar{\eps} \big( e_{\xi \ot b, \eta \ot c}(\zeta \ot 1) \big) =\bar{\eps} \left( (\xi \ot b)\langle \eta \ot c, \zeta \ot 1 \rangle \right) \\
&=& \bar{\eps} \big( \xi \ot bc^* \langle \eta , \zeta \rangle \big)=\xi \eps(bc^*)\langle \eta , \zeta \rangle = e_{\xi \eps(b c^*),\eta } (\zeta).\end{eqnarray*}
\end{proof}

The map $\hat{\eps}$ may be thought of as a kind of conditional expectation, induced by the conditional expectation $\eps$. The analogy becomes more clear, if one notices that the linear extension of the formula:
\[\langle \xi_1 \ot b_1 , \xi_2 \ot b_2 \rangle_{\eps} =\eps \left( \langle \xi_1 \ot b_1 , \xi_2 \ot b_2 \rangle \right)=\eps \left( b_1^* \langle \xi_1 , \xi_2 \rangle b_2 \right), \;\; \xi_1, \xi_2 \in E, b_1, b_2 \in B,\]
defines an $A$-linear inner product on the vector space $E \ot B$ giving rise to a Hilbert $A$-module $(E \ot B)_{\eps}$. The submodule $(E \ot 1)_{\eps}$ may be identified with $E$ and, with the corresponding identifications,  $\bar{\eps}$ projects onto it.

The construction of the map $\hat{\eps}$ may be extended in a natural way for maps between different Hilbert $A$-modules $E$, $F$ and their extensions $E \ot B$, $F \ot B$. This leads to the main theorem of this section.

\begin{tw} \label{liftable}
The quotient $\Tt_E \to \Oo_E$ is ucp liftable provided there exists a conditional expectation $\eps : F_E \to A$.
\end{tw}
\begin{proof}
We will freely use the preceding results, with $B = F_E$. Following Pimsner we denote the Hilbert bimodule $E \ot F_E$ by $E_{\infty}$. Recall first (\cite{pi}) that for each $n \in \N$ $(E \ot F_E)^{n} \cong E^n \ot F_E$, and so also $\Gamma_{E_{\infty}} \cong \Gamma_E \ot F_E$. As the existence of the lift $\Oo_E \cong \Oo_{E_{\infty}} \to \Tt_{E_{\infty}}$ has been discussed before, it is enough to check that the map $\hat{\eps} : \Ll_{F_E} (\Gamma_{E_{\infty}}) \to \Ll_A (\GE)$ acts in a proper way on the operators of the form $t_{\mu \ot b} t_{\nu \ot c}^*$, where $\mu, \nu \in \bigcup_n E^n$ and $b,c \in \Kk (E^i)$ ($i \in \N$). For $j > i$
$$\tilde{\eps} \big(e_{\mu \ot b ,\nu \ot c}  \ot \id_{E_{\infty}^j} \big) =\tilde{\eps} \left( e_{\mu ,\nu} \ot bc^* \ot \id_{E_{\infty}^{j-i}}\right) = e_{\mu ,\nu} \ot bc^* \ot \id_{E^{j-i}},$$
where (abusing the notation)  $\tilde{\eps}: \Ll_{F_E} (E^{|\mu|+j} \ot F_E , E^{|\nu|+j}\ot F_E  ) \to \Ll_A (E^{|\mu|+j} , E^{|\nu|+j})$ is a map constructed from $\eps$ as $\hat{\eps}$ was. Note that the first equality in the formula above follows from the way in which $(E \ot F_E)^{n}$ and $E^{n} \ot F_E$ are identified. Let $\pi_i$ denotes the (nonunital) representation of $\Kk(E^i)$ on $\GE$, acting only on the first $i$ components of a vector in each of the $E^j$, for $j \geq i$. Using the formula $t_{\mu}t_{\nu}^*= \sum_j e_{\mu ,\nu} \ot \id_{E_{\infty}^j}$ we find
$$\hat{\eps} \big( t_{\mu \ot b} t_{\nu \ot c}^* \big) - t_{\mu} \pi_i (bc^*) t_{\nu}^* \in \Kk (\GE),$$
which ends the proof.
\end{proof}

In general the question of the existence of the conditional expectation $\eps: F_E \to A$ is difficult. As $F_E$ is defined as the inductive limit algebra, it would be natural to construct such an expectation also inductively. If $E$ is finitely generated, each of the algebras $\Kk(E^k)$ ($k \in \N$) may be embedded in a finite matrix algebra over $A$. The embeddings are however not unique, and there is no guarantee that the standard averaging with respect to the trace would yield in a limit a conditional expectation. Even if $E$ is a finitely generated free module (i.e. $E= l_2^n \ot A$ for some $n \in \N$), so that $\Kk(E^k) \cong M_{n^k} (A)$, one has to assure that the tracial averages are compatible with the embeddings $\Kk(E^n) \hookrightarrow  \Kk(E^{n+1})$ (recall the latter are indirectly given via the left action of $A$ on $E$). Below we present a class of $C^*$-correspondences, giving rise both to crossed products and Cuntz-algebras, in which it is possible.

\medskip
\noindent
\textbf{Example.}
Let $n \in \N$, $\alpha_1, \ldots, \alpha_n \in Aut(A)$, and let $U \in M_n(A)$ be unitary. Let $E= l_2^n \ot A$ be the standard right $A$-Hilbert module and define the left action of $A$ by the formula ($ a \in A$)
\[ \phi(a) = U^* \textrm{diag} [\alpha_1(a),\ldots, \alpha_n(a)] U \in M_n(A) \cong \Kk(E).\]
In this way $E$ becomes a $C^*$-correspondence. It is easy to see that the algebra $\Oo_E$ is at the same time a certain generalisation of Cuntz algebra and the crossed product construction. When $n=1$ and $U=1$, $\Oo_E = A \rtimes_{\alpha} \Z$. If $U=1$ and $\alpha_1 = \cdots = \alpha_n = \id_A$, then $\Oo_E = A \otimes \Oo_n$ (recall that $\Oo_n$ is nuclear). When $U=1$, but $n\in \N$ is arbitrary, the algebra $\Oo_E$ is generated by Cuntz isometries $S_1, \ldots, S_n$ and a copy of $A$, with the commutation relations: $S_i a = \alpha_i(a) S_i$, valid for all $i\in \{1, \ldots, n\}$, $a \in A$. In the most general case additional twisting by $U$ is introduced.

To describe the afore-mentioned conditional expectation we need to introduce some more notation. Let $\wt{\alpha}: A \to M_n(A)$ be defined by
\[\wt{\alpha} (a) =  \textrm{diag} [\alpha_1(a),\ldots, \alpha_n(a)], \;\;\;a \in A \]
and introduce also its `inverse' $\hat{\alpha}: M_n(A) \to M_n(A)$:
\[ \hat{\alpha} ([a_{ij}]_{i,j=1}^n) = [\delta_{ij}\alpha_i^{-1} (a_{ii})]_{i,j=1}^n, \;\;\;   [a_{ij}]_{i,j=1}^n \in M_n(A)\]
($\delta$ is the standard Kronecker symbol).

Denote by $\Ex:M_n(A) \to A$ the normalised trace map:
\[ \Ex( [a_{ij}]_{i,j=1}^n ) = \frac{1}{n} \sum_{i=1}^n a_{ii},\;\;\; [a_{ij}]_{i,j=1}^n \in M_n (A).\]
Recall also the standard notation for inner automorphisms: for a unitary $V$ in a $C^*$-algebra $B$ the map $\Ad V: B \to B$ is given by $\Ad V (\cdot) = V^*  \cdot V$. In this notation the left action $\phi: A \to \Kk(E)$ may be written as
\[ \phi = \Ad U \circ \wt{\alpha} \]
and generally (for $k \geq 2$)  the embeddings $\phi_{k} : A \hookrightarrow \Kk(E^{k}) \cong M_{n^k}(A)$ are given inductively  by
\[ \phi_{k} = \Ad (I_{M_{n^k}} \ot U) \circ (I_{M_{n^k}} \ot \wt{\alpha}) \circ \phi_{k-1}
\]
(with $\phi_1=\phi$).
It is easy now to see that the required conditional expectation arises as the inductive limit of the maps $\Ex_k: M_{n^k} (A) \to A$, consecutively `inverting' the embedding procedure and applying the tracial conditional expectation. They are defined inductively  by
\[\Ex_1 = \Ex \circ \hat{\alpha} \circ \Ad U^*, \]
and for all $k \in \N$
\[ \Ex_{k+1} = \Ex_k \circ (I_{M_{n^k}} \ot \Ex )\circ (I_{M_{n^k}} \ot \hat{\alpha}) \circ \Ad (I_{M_{n^k}} \ot U^*).\]

\section{Approximation Properties}         \label{appprop}

There are at least four basic approximation properties for a $C^*$-algebra $A$ considered in the literature, all closely related to the minimal tensor product (denoted here by $\otimes$). These are in increasing order of generality:
\newcounter{Lcount}  
\begin{list}
{\arabic{Lcount}.}
{\usecounter{Lcount}}
\item  Nuclearity, which is equivalent to the CPAP (completely positive approximation property): there exists a net of completely positive contractions $\varphi_{\lambda}:A \to M_{n_{\lambda}}$ and $\psi_{\lambda} : M_{n_{\lambda}} \to A$ such that $\psi_{\lambda} \circ \varphi_{\lambda}(x) \to x$ for all $x \in A$.
\item  The CBAP (completely bounded approximation property): there exists a net $(\phi_{\lambda}:A \to A)$ of finite rank maps such that $\phi_{\lambda}(x) \to x$ for all $x \in A$ and $\sup_{\lambda}  \| \phi_{\lambda} \|_{cb} < \infty$. The smallest possible such supremum is the Haagerup constant $\Lambda(A)$ of $A$.
\item The strong OAP (strong operator approximation property): there exists a net $(\phi_{\lambda}:A \to A)$ of finite rank maps such that $(\phi_{\lambda} \ot \textup{id})(x) \to x$ for all $x \in A \ot B(\ell^2(\N))$.
\item Exactness, which is equivalent to nuclear embeddability: for every faithful representation $A \to B(H)$ there exists a net of completely positive contractions $\varphi_{\lambda}:A \to M_{n_{\lambda}}$ and $\psi_{\lambda} : M_{n_{\lambda}} \to B(H)$ such that $\psi_{\lambda} \circ \varphi_{\lambda}(x) \to x$ for all $x \in A$.\end{list}
Finally there is
\begin{list}{\arabic{Lcount}.}{\usecounter{Lcount}}\setcounter{Lcount}{4}
\item The OAP (operator approximation property), which neither implies nor follows from exactness: there exists a net $(\phi_{\lambda}:A \to A)$ of finite rank maps such that $(\phi_{\lambda} \ot \textup{id})(x) \to x$ for all $x \in A \ot \Kk(\ell^2(\N))$\end{list}
Note however that it is known that a $C^*$-algebra has strong OAP if and only if it is exact and has OAP (\cite{HaagerupKraus} Thm.2.2).

The last three properties allow a characterisation via Fubini products. Given a closed subspace $S$ of a $C^*$-algebra $B$ recall that the Fubini product is defined by
$$F(A,S)=F(A,S,A \ot B)=\{ x \in A \ot B \mid \forall \omega \in A^* : (\omega \ot \textup{id})(x) \in S \}.$$
$A$ has the (strong) OAP iff $F(A,S,A \ot \Kk(\ell^2))= A \otimes S$ ($F(A,S,A \ot B(\ell^2))= A \otimes S$) for all closed subspaces $S$ of $\Kk(\ell^2)$ ($B(\ell^2)$) and $A$ is exact iff $F(A,J,A \ot E)=A \ot J$ for all $C^*$-algebras $E$ and ideals $J \subs E$. The following lemma shows that to prove an approximation property one may replace approximating finite-rank factorisations with the factorisations via algebras possessing the property in question (a similar result concerning exactness appears in \cite{Dykema}).

\begin{lem} \label{app}
Suppose there exists an approximating net $(\varphi_i : A \to C_i, \psi_i : C_i \to A)$ i.e.\ $\psi_i \circ \varphi_i (a)\to a$ for all $a \in A$, where $\varphi_i$ and $\psi_i$ are contractive and completely positive. If for any of the five approximation properties all $C_i$ have this property then so does $A$, except in case of the CBAP, where $A$ has the OAP if all $C_i$ have the CBAP and $A$ has CBAP if $\sup_i \Lambda (C_i) < \infty$.
\end{lem}
\begin{proof}
We only show the claim for exactness, the other cases are obvious. In this situation we need to show that $F(A,J,A \ot E)=A \ot J$ for all $C^*$-algebras $E$ and ideals $J \subs E$. Now $(\varphi_i \ot \textup{id}_E ) (F(A,J)) \subs F(C_i , J) $ since $(\omega \ot \textup{id}_E) (\varphi_i \ot \textup{id}_E)(x) = (( \omega \circ \varphi_i ) \ot \textup{id}_E )(x) \in J$ for $x \in F(A,J)$ and $\omega \in C_i^*$. By exactness of $C_i$,  $(\varphi_i \ot \textup{id}_E)(x) \in C_i \ot J$. The latter is mapped into $A \ot J$ under $\psi_i \ot \textup{id}_E$. Since $(\psi_i \circ \varphi_i) \otimes \textup{id}_E (x) \to x$ for all $x \in A \otimes E$ the result follows.
\end{proof}

The two following propositions are undoubtedly well known, but we were not able to locate a convenient reference.

\begin{prop}
 Let $A$ be a $C^*$-algebra with a (closed, two-sided) ideal $J$ and let (*) be one of the properties 1.--5. above. If $A$ has (*), then $J$ also has (*).
\end{prop}
\begin{proof}
It is enough to use Lemma \ref{app} and the approximating maps given by the twosided multiplication by a contractive approximate unit in $J$.
\end{proof}

\begin{prop} \label{compact}
Let $A$ be a $C^*$-algebra, $Y$ a Hilbert $C^*$-module over $A$ (not necessarily full) and let (*) be one of the properties 1.--5.\ above. If $A$ has (*), then $\Kk(Y)$ also has (*).
\end{prop}
\begin{proof} The (completely) contractive maps giving approximate factorisations of a Hilbert module $Y$ via the column modules $A^{\oplus n}$ ($ n \in \N$) constructed in \cite{Ble} induce (by twosided multiplications) the completely positive and contractive factorisations of $\Kk(Y)$ via $\Kk(A^{\oplus n}) \cong M_n(A)$ (see \cite{ds} - note however it is not necessary to assume that $Y$ is full or countably generated). As it is clear that the approximation properties of $M_n(A)$ are the same as those of $A$, Lemma \ref{app} ends the proof.
\end{proof}

For the sake of completeness we include one more fact (which is not actually used anywhere else in the paper).

\begin{prop} \label{extens}
Let (*) be one of the properties 1.--5. above and let $ 0 \to J \to A \stackrel{q}{\to} B \to 0 $ be a short exact sequence of $C^*$-algebras. Assume that $q$ is locally liftable; that is for each finite dimensional operator system $X \subset B$ there exists a ucp map $\varphi: X \to A$ such that $q \circ \varphi = \textrm{id}_X$. Then if $J$ and $B$ have (*), then also $A $ has (*).
\end{prop}
\begin{proof}
In case of nuclearity and exactness the conclusion can be reached along the lines of Proposition 2.1 of \cite{Joachim}, where the existence of quasicentral approximate units in $J$ is exploited to construct relevant approximating nets. For the CBAP similar arguments as in [DyS, Thm.1] can be used. To apply the construction in case of the OAP we need local liftability of $q \otimes \textrm{id}_{\Kk}$ assuming local liftability of $q$. This follows immediately from the Effros-Haagerup lifting theorem ([EfH, 3.2]): all we need to show is that $ 0 \to J\otimes \Kk \ot C  \to A \otimes \Kk \ot C \to B\otimes \Kk \ot C \to 0 $ is exact for every $C^*$-algebra $C$. But $ 0 \to J \ot C  \to A \otimes  C \to B \otimes   C \to 0$ is exact and $\Kk$ is an exact $C^*$-algebra. Finally the case of strong OAP follows since strong OAP is equivalent to OAP together with exactness, as pointed out  after the introduction of approximation properties.
\end{proof}

Combining the above with the results of previous sections we obtain the first of two main theorems of the paper:

\begin{tw}
Let (*) be one of the properties 1.--5. above and let  $X$ be a finitely generated Hilbert bimodule over a unital $C^*$-algebra $A$. If $A$ has (*) then the Cuntz-Pimsner algebra $\Oo_X$ has (*).
\end{tw}

As was mentioned in Section 3 Pimsner proved in \cite{pi} that whenever $E$ is a $C^*$-correspondence, then $\Oo_E \cong \Oo_{E_{\infty}}$, where $E_{\infty}$ is a $C^*$-bimodule over $F_E$ (see also \cite{AbAch}). As $F_E$ is an inductive limit of the $C^*$-algebras $\Kk(E^n)$, Proposition \ref{compact} and permanence of approximation properties with respect to taking inductive limits imply the following corollary:

\begin{cor} \label{correspond}
Let (*) be one of the properties 1.--5.\ above. If $E$ is a finitely generated $C^*$-correspondence over a unital $C^*$-algebra $A$ and $A$ has (*), then the Cuntz-Pimsner algebra $\Oo_E$ has (*).\end{cor}

If there exists a lifting our methods can be also applied to the case of $\TE$.

\begin{tw} Let (*) be one of the properties 1.--5.\ above and let  $E$ be a finitely generated $C^*$-correspondence over a unital $C^*$-algebra $A$. If $A$ has (*) and the short exact sequence $ 0 \to \Kk (\Gamma_E) \to \TE \to \Oe \to 0 $ is ucp liftable, then both the Toeplitz-Pimsner algebra $\Tt_E$ and the Cuntz-Pimsner algebra $\Oo_E$ have (*).
\end{tw}

Recall that the existence of the required lift is guaranteed if there exists a conditional expectation $\eps : F_E \to A$ (Theorem \ref{liftable}), or if $A$ is nuclear, separable and $E$ is countably generated (then $\Oo_E$ is separable, by Corollary \ref{correspond} it is also nuclear, and one can apply the Choi-Effros lifting theorem (\cite{ChE})). In fact the implication $A$-nuclear $\Longrightarrow$ $\TE$-nuclear holds even without any separability assumptions, as the existence of local lifts is sufficient to repeat the argument.

\end{document}